\documentclass[a4paper,11pt]{article}
\usepackage{amsmath, amsfonts}
\usepackage[amsmath,thmmarks]{ntheorem}
\usepackage{enumitem,tikz}
\usetikzlibrary{arrows}
\title{Vahlen groups defined over commutative rings}
\author{Justin M\textsuperscript{c}Inroy \\ \footnotesize{Heilbronn Institute for Mathematical Research, School of Mathematics,} \\ \footnotesize{University of Bristol, University Walk, Bristol, BS8 1TW, UK} \\ \footnotesize{Tel: +44 117 954 5661, Fax: +44 117 331 5264,} \\ \footnotesize{email: justin.mcinroy@bristol.ac.uk}}

\newcommand\la{\langle}
\newcommand\ra{\rangle}

\newtheorem{thm}{Theorem}[section]
\newtheorem{lem}[thm]{Lemma}
\newtheorem{cor}[thm]{Corollary}
\newtheorem{prop}[thm]{Proposition}

\theorembodyfont{\normalfont}
\newtheorem{defn}[thm]{Definition}

\makeatletter
\newtheoremstyle{MyNonumberplain}%
  {\item[\theorem@headerfont\hskip\labelsep ##1\theorem@separator]}%
  {\item[\theorem@headerfont\hskip\labelsep ##3\theorem@separator]}
\makeatother
\theoremstyle{MyNonumberplain}

\theoremheaderfont{\normalfont\itshape}
\theorembodyfont{\normalfont}
\theoremseparator{.}
\theoremsymbol{\ensuremath{\square}}
\newtheorem{pf}{Proof}
\qedsymbol{\ensuremath{\square}}

\theoremheaderfont{\normalfont\bfseries}
\theorembodyfont{\itshape}
\theoremseparator{}
\theoremsymbol{}

\newcommand{\T}{\mathcal{T}}
\renewcommand{\phi}{\varphi}

\begin{document}

\maketitle

\begin{abstract}
Elements of a Vahlen group are $2 \times 2$ matrices with entries in a Clifford algebra satisfying some conditions.  Traditionally they have come in both ordinary and paravector type and have been defined (over Clifford algebras) over the real or complex numbers.  We extend the definition of both types to be over a commutative ring with an arbitrary quadratic form.  We show that they are indeed groups and identify in each case the group as the pin group, spin group, or another subgroup of the Clifford group.  Under some mild conditions, for both types we show the equivalence of our definition with a suitably generalised version of the two standard definitions.
\end{abstract}

\section{Introduction}

It is well-known that M\"obius transformations can be represented by $2 \times 2$ matrices.  In 1902, Vahlen extended this to higher dimensions by defining Vahlen groups \cite{vahlen}.  He described a group whose elements were $2 \times 2$ matrices with entries in a real Clifford algebra and satisfying some conditions.  Using the upper half-plane model of hyperbolic space, one can then show that Vahlen's group acts on this and is a cover of the group of orientation preserving isometries.

In fact, what Vahlen introduced and was later re-presented by Ahlfors in \cite{ahlfors2} is what we call paravector Vahlen groups.   Ordinary Vahlen groups were introduced by Maa\ss{} \cite{maass} (see also \cite{ahlfors1}).  The paravector Vahlen groups are a cover of the orientation preserving M\"obius transformations, while the ordinary Vahlen groups are a cover of the full group of M\"obius transformations \cite[Chapter 19]{lounesto}.  The ordinary Vahlen groups defined there correspond to the case when the quadratic form is definite.  The paravector groups correspond to a form of signature $(+1,-n)$.  If we take an $n$-dimensional vector space $V$ with a negative definite quadratic form $q$, then a \emph{paravector} is an element of $\mathbb{R} \oplus V$.  Using paravectors, there is a natural extension of $q$ to a quadratic form of signature $(+1,-n)$.  In particular, taking $\dim V = 3$, we get Minkowski space-time.  We note that for both types there are two common but slightly different equivalent definitions which may be given.

One can also define (paravector) Vahlen groups for indefinite forms.  These are more difficult and were tackled by Cnops in \cite{cnops}.  Complex Vahlen groups can also be defined.  In all real and complex cases, the group is isomorphic to the pin group, spin group, or a related subgroup of the Clifford group.  A survey of such results can be found, for example, in \cite[Chapter 19]{lounesto} or \cite[Chapter 18]{porteous}.

In 1987, Elstrodt, Grunewald and Mennicke defined (special paravector) Vahlen groups for fields of odd characteristic and for an arbitrary quadratic form \cite{vahlenpodd}.  They used a different definition to before, but showed that the group they obtained was isomorphic to the expected spin group.  With an assumption on the quadratic form, they are able to show that their definition is equivalent to the two previous definitions.

In this paper, we generalise Elstrodt, Grunewald and Mennicke's definition \cite{vahlenpodd} to define both ordinary and paravector Vahlen groups in their most general setting.  That is, over Clifford algebras $Cl(N,q)$, where $N$ is a module defined over a commutative ring $R$ and $q$ is an arbitrary (not necessarily non-degenerate) quadratic form.  We define the general ordinary Vahlen group $GV(N,q)$ to be the set of matrices
\[
 \begin{pmatrix} \alpha & \beta \\ \gamma & \delta \end{pmatrix}  \in {\rm MAT}_2(Cl(N,q))
\]
subject to some relations between the entries (see Definition \ref{vahlendef} for details).  The special ordinary Vahlen group $SV(N,q)$ is defined by requiring the pseudo-determinant to be $1$.  The general and special paravector Vahlen groups, $GPV(L,q)$ and $SPV(L,q)$, are defined similarly.  We note that in the literature, there is no agreement on notation.  In particular, authors often do not make explicit whether they consider ordinary or paravector Vahlen groups and sometimes they also only consider the special Vahlen group.  We endeavour to be clear in this paper.

In Theorems \ref{vahleniso} and \ref{paravahleniso}, we show
\begin{align*}
GV(N,q) &\cong \{ u \in C(M,q) : u\bar{u} \in R^* \}\\
SV(N,q) &\cong Pin(M,q)\\
GPV(L,q) &\cong GV(L \perp \la z \ra,q)_0\\
SPV(L,q) &\cong SV(L \perp \la z \ra,q)_0 \cong Spin(M,q)
\end{align*}
where $N$ and $L$ are free modules, and $M$ is decomposed as $M = H \perp N$ with $H$ a 2-dimensional hyperbolic module for ordinary Vahlen groups and decomposed as $M = H \perp \la z \ra \perp L$ with $H$ is a 2-dimensional hyperbolic module and $q(z) = -1$ for paravector Vahlen groups.

We also give two further definitions for the ordinary and paravector Vahlen groups and show their equivalence to our definition.  These definitions are generalised versions of the usual ones first given by Vahlen \cite{vahlen}, Maa\ss{} \cite{maass} and Cnops \cite{cnops}, and involve defining a set in which the entries of the matrix must lie.  These two definitions are shorter, but our first definition is easier to work with and is crucial for proving the isomorphisms above.

In the original case, the form was definite (or in the paravector case, built from a definite form) and so the entries of the matrix are invertible.  However, in general this is not true.  Indeed, for the indefinite real case Maks gave an example of a Vahlen matrix where none of the entries are invertible \cite[p. 41]{maks}.  We introduce a monoid $\T(N,q)$ (respectively $\mathcal{PT}(L,q)$ for the paravector case) in which the entries lie.  Under the assumption that $R$ is an integral domain and that the monoid $\T(N,q)$ (respectively $\mathcal{PT}(L,q)$) is closed under the transpose map, we show that all three definitions given are equivalent.  This is Theorem \ref{vahlenequiv} for ordinary Vahlen groups and Theorem \ref{paravahlenequiv} for the paravector Vahlen groups.  We note that this generalises a result of Elstrodt, Grunewald and Mennicke on a special class of quadratic forms over a field of odd characteristic.  Our proof covers their case and we believe it is also less computational and more elegant.

Our original motivation for studying Vahlen groups came from a conjecture of Sidki about the finiteness of a family of presentations for group $y(m,n)$.  The author together with Shpectorov and Sidki conjecture that these groups are actually orthogonal groups of dimension $m+1$ over a field $GF(2^{\frac{n-1}{2}})$.  In \cite{sidki}, the author together with Shpectorov generalise the original presentation to one for a group $y(m)$, where the groups $y(m,n)$ can be recovered as quotients of $y(m)$.  They give two homomorphisms from $y(m)$, one into the orthogonal group of dimension $m+1$ over a ring of Laurent polynomials $\mathbb{F}_2[t,t^{-1}]$ and one into the corresponding Clifford algebra.  This Clifford algebra satisfies the conditions of Theorem \ref{vahleniso} and hence we may map $y(m)$ into the, now suitably generalised, Vahlen group.  However, it can be shown that it is a proper subgroup, so further investigation, possibly using congruence subgroups, is needed.  We note that Halkjaer \cite{Claus} used Vahlen groups over fields of characteristic $p$ to solve the analogous problem for $y_p(m,n)$.

In Section \ref{sec:background}, we introduce Clifford algebras and briefly describe their properties.  We also define some related groups: the Clifford group, pin group and spin groups.  Ordinary and paravector Vahlen groups are defined in Sections \ref{sec:vahlen} and \ref{sec:paravahlen}, respectively.  We show that they are indeed groups and in all cases identify them.  Sections \ref{sec:vahlenequiv} and \ref{sec:paravahlenequiv} are devoted to showing the equivalence of different definitions for ordinary and paravector Vahlen groups, respectively.

\section{Background on Clifford algebras}\label{sec:background}

We will review the definitions and basic properties for Clifford algebras and the pin and spin groups.  For more details refer, for example, to \cite[Chapter 7]{omeara}.

Let $M$ be a module over a commutative ring $R$ and let $q:M \to R$ be a quadratic form with associated bilinear form $(\cdot, \cdot )$.  Given a bilinear form $(\cdot,\cdot): M \times M \to R$, by fixing one or other of it's variables $m \in M$ we get a map $M \to Hom(M,R)$.  We say that $(\cdot,\cdot)$ is \emph{non-degenerate} if for all $m \in M$ the map obtained is an isomorphism.  We say that $q$ is \emph{non-degenerate} if the bilinear form $(\cdot,\cdot)$ is non-degenerate.  

Let $M$ be a module over a commutative ring with a quadratic form $q:M \to R$ (not necessarily non-degenerate).  If $A$ is a unital associative algebra over $R$ and $j:M \to A$ is a linear map such that
\[
j(m)^2 = q(m).1_A
\]
for all $m \in M$, then we say that $(A,j)$ is a \emph{compatible} with the pair $(M,q)$.  The \emph{Clifford algebra} $(Cl,i)$ is defined to be the universal such pair.  That is, given any compatible pair $(A,j)$, there exists a unique algebra homomorphism $f:Cl \to A$ such that the following diagram commutes:
\begin{center}
\begin{tikzpicture}[scale=4]
\node (1) at (0,1) {$M$};
\node (2) at (1,1) {$Cl(M,q)$};
\node (3) at (0.5,0.13) {$A$};
\path[->, >=angle 90]

(1) edge node[anchor = south]{$i$} (2)
(1) edge node[anchor = north east]{$j$} (3)
(2) edge node[anchor = north west]{$f$} (3);
\end{tikzpicture}
\end{center}
We will often just write $Cl(M,q)$, or $Cl$ for the pair $(Cl,i)$.  We will also say that $j:M \to A$ \emph{extends} to a map $f:Cl(M,q) \to A$.

We may construct the Clifford algebra $Cl(M,q)$ as the quotient algebra $T(M)/I$ of the tensor algebra $T(M)$ by the ideal $I$ generated by the relations
\[
m^2 = q(m)
\]
for all $m \in M$.  If we apply the above relation to $m+n$, then we obtain
\[
mn + nm = (m,n)
\]
for all $m,n \in M$.  Viewed in this way, $M$ and $R$ embed naturally in the algebra $Cl$.  We will abuse notation and say that $R$ and $M$ lie in $Cl$.  

Since the tensor algebra has an $\mathbb{N}$-grading given by the rank of tensors, the quotient $Cl$ inherits a $\mathbb{Z}_2$-grading.  That is, $c \in Cl$ is in the even part, notated by $Cl_0$, if $c$ is the sum of tensors of even rank and in the odd part, $Cl_1$, if it is the sum of tensors of odd rank.  This has the property that $Cl = Cl_0 \oplus Cl_1$ with the $Cl_i$ being submodules, $R \subset Cl_0$ and $Cl_iCl_j \subseteq Cl_{i+j}$ for all $i,j \in \mathbb{Z}_2$.  We note that this implies that $Cl_0$ is a subalgebra of $Cl$.

The algebra has some natural automorphisms and anti-automorphisms.  The map $-1_M:M \to M$ can be extended to an involutory automorphism on $Cl$ which we write as $x'$, for $x \in Cl$; this is called the \emph{grade automorphism}.  If $x$ has graded decomposition $x = x_0 + x_1 \in Cl_0 \oplus Cl_1$, then $x' = (x_0+x_1)' = x_0 - x_1$.  There is also a \emph{transpose} map from $Cl$ to $Cl^{op} \cong Cl$ given by reversing the order of elements in products, i.e., it is defined by extending linearly the map $(m_1 \dots m_k)^t = m_k \dots m_1$, where $m_i \in M$. Clearly, the transpose map satisfies $(xy)^t = y^t x^t$, for $x,y \in Cl$, making it an involutory anti-automorphism of $Cl$.  Combining these two maps, we get the \emph{Clifford conjugation} map $\overline{x} = (x^t)' = (x')^t$ which is also an involutory anti-automorphism of $Cl$.  Note that in characteristic two, Clifford conjugation is just the transpose map.

Associated with the grade automorphism we can define a subalgebra of $Cl$ which is a generalisation of the centre.  We call the subalgebra
\[
TCen(Cl):=\{ x \in Cl : xm = mx' \quad \forall m \in M \}
\]
the \emph{twisted centre} of $Cl$.  We note that elements of the twisted centre commute with $a \in Cl_0$, but `twisted commute' with elements of $Cl_1$.  Also, $TCen(Cl)_0 \subseteq Cen(Cl)$.

Let $N$ denote the map $N(x) = x\bar{x}$ for $x \in Cl$.  Suppose that $a,b \in Cl$ such that $N(a),N(b) \in Cen(Cl)$.  Then, $N(ab) = ab \bar{b}\bar{a} = a\bar{a} b \bar{b} = N(a)N(b)$.  So, if $N$ is restricted to $\{x \in Cl : N(x) \in Cen(Cl) \}$, then $N$ is a homomorphism which we will call the \emph{norm}; this is an extension of the usual definition.  One can also see by a simple calculation that $N(x)$ is even.  We note that if $x$ invertible, then $\bar{x}$ and so $N(x)$ is also invertible.  Then,
\[
1 = u^{-1}u = \bar{u}N(u)^{-1}u = N(\bar{u})N(u)^{-1}
\]
and so $N(\bar{u}) = N(u)$.  However, this property is not known to hold in general for non-invertible elements.

There are several groups associated with the Clifford algebra $Cl(M,q)$.  The set of all units of $Cl$ form a group, notated $Cl(M,q)^\times$.  For $u \in Cl^\times$, define $\pi(u)$ as the map $x \mapsto u'xu^{-1}$.  Then, we define the \emph{Clifford group} by
\[
C(M,q):= \{ u \in Cl(M,q)^\times : \pi(u)m \in M \quad \forall m \in M\}
\]
When restricted to $C(M,q)$, $\pi$ is a natural homomorphism to $GO(M,q)$.  If $n \in M$ such that $q(n) \in R^*$, where $R^*$ is the set of units of $R$, then $n$ is in $Cl^\times$ with $n^{-1} = n/q(n)$.  Since $mn+nm = (m,n)$, we have that
\[
\pi(n)m = n'mn^{-1} = (mn - (n,m))n^{-1}  = m - \frac{(n,m)}{q(n)}n.
\]
Hence, $n \in C(M,q)$ and $\pi(n)$ is an orthogonal reflection $r_n$ (without the use of $n'$ in the definition it would be $-r_n$).  In particular, if $R$ is a field, then $\pi$ is onto $GO(M,q)$.

\begin{lem}\label{piiff}
Let $u \in C(M,q)$, $x \in Cl$.  Then, $\pi(u)x \in M$ if and only if $x \in M$.
\end{lem}
\begin{pf}
Apply $\pi(u)^{-1}$ to get the desired result.
\end{pf}

The kernel of $\pi$ is the set $TCen(Cl)^\times$ of elements of the twisted centre $TCen$ which are invertible.  Hence, we see that $TCen^\times$ is a normal subgroup of $C(M,q)$.  We note that if $q$ is non-degenerate, then $\ker \pi = TCen(Cl)^\times = R^\times$ (see, for example, \cite[Section 2.4]{bass}).

The Clifford group is closed under the grade automorphism, transpose map and conjugation.  In fact, we have the following:

\begin{lem}\label{C(M,q)}{\rm \cite[Proposition 3.2.1]{bass}}
If $u \in C(M,q)$, then $\bar{u} \in C(M,q)$.  Moreover, $N(u) = u \bar{u} \in \ker \pi$ and so $\pi(\bar{u})$ is in the preimage $\pi(u)^{-1}$.
\end{lem}

We note however, that in general for $u \in C(M,q)$, the norm $N(u) \in Cen_0^\times$ and will not necessarily lie in $R$.

We have the following simple lemma.

\begin{lem}\label{piequiv}
If $u$ is an invertible element of $Cl$, then $u^{-1} = \bar{u}N(u)^{-1}$.  Moreover, for $x \in Cl$,
\[
\pi(u)x = u'xu^{-1} = u'x\bar{u}N(u)^{-1} = (ux'u^t)'N(u)^{-1}.
\]
\qed
\end{lem}

In particular, if $u$ is invertible, then when checking that $\pi(u)m \in M$ we may instead check that $umu^tN(u)^{-1} \in M$ for all $m \in M$.  In particular, if $N(u) \in R^*$ and $R$ is an integral domain, we need only check $umu^t$.  Also, the above gives us the following.

Since $N$ restricted to $C(M,q)$ is a homomorphism, we may now define the \emph{pin group} as the kernel of this homomorphism.
\[
Pin(M,q):=\{ u \in C(M,q) : N(u) = u \bar{u}=1 \}
\]
Hence, we have $Pin(M,q) \unlhd C(M,q)$ and $C(M,q)/Pin(M,q)$ is isomorphic to a subgroup of $TCen^\times$.  We note that when $q$ is non-degenerate, then $TCen^\times$ can be replaced with $R^*$ \cite[Proposition 3.2.1]{bass}.

We may now use the grading of $Cl$ to define the remaining groups.  We define the \emph{even Clifford group} to be the even graded part of $C(M,q)$.
\[
C^+(M,q) :=C(M,q)_0
\]
Finally, we define the \emph{spin group} to be the even graded part of the pin group.  It is also the kernel of $N$ restricted to $C^+(M,q)$.
\begin{align*}
Spin(M,q) :=& Pin(M,q)_0 \\
=& \ker N|_{C^+(M,q)}
\end{align*}

\section{Vahlen groups}\label{sec:vahlen}

In this section, we will define Vahlen groups with a (possibly degenerate) quadratic form over a commutative ring.  This generalises a definition of Elstrodt, Grunewald and Mennicke in \cite{vahlen}.  We will show that a Vahlen group is indeed a group and identify it as a subgroup of the Clifford group.

\begin{defn}\label{vahlendef}
Let $N$ be a module over a commutative ring $R$ with a quadratic form $q$.  Define the \emph{general Vahlen group} $GV(N,q)$ to be the set of matrices
\[
\begin{pmatrix} \alpha & \beta \\ \gamma & \delta \end{pmatrix}  \in {\rm MAT}_2(Cl(N,q))
\]
such that the following conditions hold.
\begin{enumerate}[label ={\rm (\roman*)}]
	\item $\alpha \delta^t - \beta \gamma^t \in R^*$
	\item $\alpha \beta^t = \beta \alpha^t$ and $\gamma \delta^t = \delta \gamma^t$
	\item $\alpha \bar{\alpha}, \beta\bar{\beta}, \gamma \bar{\gamma}, \delta \bar{\delta} \in R$
	\item $\alpha\bar{\gamma}, \beta\bar{\delta} \in N \cup \{0\}$
	\item for all $x \in N$, $\alpha x \bar{\beta} + \beta \bar{x} \bar{\alpha}, \gamma x \bar{\delta} + \delta \bar{x} \bar{\gamma} \in R$
	\item for all $x \in N$, $\alpha x \bar{\delta} + \beta \bar{x} \bar{\gamma} \in N \cup \{0\}$
\end{enumerate}
We call $\alpha \delta^t - \beta \gamma^t$ the \emph{pseudo-determinant} and define the \emph{special Vahlen group} $SV(N,q)$ to be the subset of $GV(N,q)$ with pseudo-determinant $1$.
\end{defn}

We will establish the following:

\begin{prop}\label{vahlengroup}
The set $GV(N,q)$ is a group with inverses given by
\[
\begin{pmatrix} \alpha & \beta \\ \gamma & \delta \end{pmatrix}^{-1} = \frac{1}{\alpha \delta^t - \beta \gamma^t} \begin{pmatrix} \delta^t & -\beta^t \\ -\gamma^t & \alpha^t \end{pmatrix}.
\]
Moreover, the pseudo-determinant is a homomorphism into $R^*$ whose kernel is $SV(N,q)$.
\end{prop}

In fact we will show more, that $GV(N,q)$ and $SV(N,q)$ are isomorphic to $NC(M,q)$ and $Pin(M,q)$, respectively, for some module $M$, where
\[
NC(M,q) :=\{ u \in C(M,q) : N(u) \in R^* \}
\]

Let $M$ be a module over $R$ with a quadratic form $q$.  Suppose that there exists a splitting $M = \la e,f \ra \perp N$ such that $\la e,f \ra$ is a hyperbolic line and $N$ is a free module.  It is well-known that, given such a splitting, the Clifford algebra is isomorphic to a matrix algebra over a smaller Clifford algebra:

\begin{lem}{\rm \cite[Lemma 7.1.9]{omeara}}
Suppose $M$ has an orthogonal splitting $M=H\perp N$, where $H$ is a free hyperbolic submodule of dimension $2k$ and $N$ is free of finite rank.  Then,
\[
Cl(M,q) \cong {\rm MAT}_{2^k}(Cl(N,q))
\]
\end{lem}

We will briefly sketch the proof for $k=1$ which is our situation.  This will also allow us to explicitly describe the isomorphism and fix some notation.  Define a map $\phi: M \to {\rm MAT}_{2}(Cl(N,q))$ by
\[
\phi(er+fs + n) = \begin{pmatrix} n & r \\ s & -n \end{pmatrix}
\]
where $r,s \in R$ and $n \in N$.  Then,
\begin{align*}
\phi(er+fs + n)^2 &= \begin{pmatrix} n & r \\ s & -n \end{pmatrix}^2 \\
&= \begin{pmatrix} n^2 + rs & 0 \\ 0 & n^2 + rs \end{pmatrix} \\
&= (n^2+rs) \begin{pmatrix} 1 & 0 \\ 0 & 1 \end{pmatrix}
\end{align*}
So, since $q(er+fs + n) = n^2+rs$, $({\rm MAT}_{2}(Cl(N,q)), \phi)$ is a compatible pair and so $\phi$ extends to a unique algebra homomorphism from $Cl(M,q)$.  By considering the dimensions of the algebras, we see that it is in fact a isomorphism.  By an abuse of notation, we will use the same notation for this giving us
\[
\phi: Cl(M,q) \to {\rm MAT}_{2}(Cl(N,q))
\]
We observe that all elements of $Cl(M,q)$ can be written in the form $ef \alpha + e \beta' + f \gamma + fe \delta'$, where $\alpha, \beta, \gamma, \delta \in Cl(N,q)$.  Since, $e$ and $f$ are both orthogonal to $N$, note that
\[
ef \alpha + e \beta' + f \gamma + fe \delta' = \alpha ef + \beta e + \gamma' f  + \delta' fe
\]
\begin{lem}\label{matiso}
The isomorphism $\phi: Cl(M,q) \to {\rm MAT}_{2}(Cl(N,q))$ is given by
\[
\phi(ef \alpha + e \beta' + f \gamma + fe \delta') = \phi(\alpha ef + \beta e + \gamma' f  + \delta' fe) = \begin{pmatrix} \alpha & \beta \\ \gamma & \delta \end{pmatrix}
\]
where $\alpha, \beta, \gamma, \delta \in Cl(N,q)$.
\end{lem}
\begin{pf}
It is an easy calculation to show that the above map is an algebra homomorphism.  Noting that $n = (ef+fe)n$ and $n' = -n$, we see that $\phi$ as a map defined from $M$ agrees with the above map restricted to $M$.  However, by the universality of Clifford algebras, the isomorphism $\phi: Cl(M,q) \to {\rm MAT}_{2}(Cl(N,q))$ is the unique extension of such a map.  Hence, $\phi$ is given by the above map.
\end{pf}

We will now use this isomorphism to write the standard Clifford automorphisms and anti-automorphisms in terms of matrices.

\begin{lem}\label{translate}
Let $\begin{pmatrix} \alpha & \beta \\ \gamma & \delta \end{pmatrix} = \phi(x)$ for some $x \in Cl(M,q)$.  Then,
\begin{align*}
\phi(x') &= \begin{pmatrix} \alpha' & -\beta' \\ -\gamma' & \delta' \end{pmatrix} \\
\phi(x^t) &= \begin{pmatrix} \bar{\delta} & \bar{\beta} \\ \bar{\gamma} & \bar{\alpha} \end{pmatrix} \\
\phi(\bar{x}) &= \begin{pmatrix} \delta^t & -\beta^t \\ -\gamma^t & \alpha^t \end{pmatrix}
\end{align*}
We will define these to be $\phi(x)'$, $\phi(x)^t$ and $\overline{\phi(x)}$, respectively.
\end{lem}
\begin{pf}
We do just the first of these calculations, the rest are similar.
\begin{align*}
\phi(x') &= \phi\left((ef \alpha + e \beta' + f \gamma + fe \delta')'\right) \\
&=\phi(ef \alpha' - e \beta - f \gamma' + fe \delta) \\
&=\begin{pmatrix} \alpha' & -\beta' \\ -\gamma' & \delta' \end{pmatrix}
\end{align*}
\end{pf}

Note that the group of invertible elements of $Cl(M,q)$ is isomorphic to $GL_2(Cl(N,q))$.

\begin{thm}\label{vahleniso}
Let $M$ be a module over $R$ with a quadratic form $q$.  Suppose that there exists a splitting $M = \la e,f \ra \perp N$ such that $\la e,f \ra$ is a hyperbolic line and $N$ is a free module.  Then,
\begin{align*}
GV(N,q) &\cong  NC(M,q)\\
SV(N,q) &\cong Pin(M,q)
\end{align*}
\end{thm}

\begin{pf}
Using Lemma \ref{translate}, we begin by considering the following.
\begin{align}
\phi(u\bar{u}) &= \begin{pmatrix} \alpha & \beta \\ \gamma & \delta \end{pmatrix} \begin{pmatrix} \delta^t & -\beta^t \\ -\gamma^t & \alpha^t \end{pmatrix} \nonumber \\
&=\begin{pmatrix} \alpha \delta^t - \beta \gamma^t & -\alpha \beta^t+ \beta \alpha^t \\ \gamma \delta^t - \delta \gamma^t & \delta\alpha^t - \gamma \beta^t \end{pmatrix} \label{ubaru}
\end{align}
If $u \in NC(M,q)$, then $u$ is invertible and $N(u) = u \bar{u} \in R^*$.  Conversely, if $N(u) \in R^*$, then by the above calculation, $\phi(u)$ is an invertible matrix and hence $u$ is also invertible.  By Lemma \ref{matiso}, $r:=N(u) \in R^*$ is equivalent to $\phi(r) = \begin{pmatrix} r & 0 \\ 0 & r \end{pmatrix}$.  This is true if and only if $\alpha \beta^t = \beta \alpha^t$, $\gamma \delta^t = \delta \gamma^t$ and $\alpha \delta^t - \beta \gamma^t \in R^*$.  Note that, if  $\alpha \delta^t - \beta \gamma^t \in R^*$, then $\alpha \delta^t - \beta \gamma^t = (\alpha \delta^t - \beta \gamma^t)^t = \delta\alpha^t - \gamma \beta^t$.  So $u$ being invertible and $N(u) \in R^*$ is equivalent to conditions (3)(i) and (ii).  Also note that the pseudo-determinant is $1$ precisely when $N(u) = u \bar{u} = 1$.

Since we may now assume that $u$ is invertible and $N(u) \in R^*$, $u^{-1} = \bar{u}/N(u)$ and to show $\pi(M) \subset M$ we need only consider $umu^t$.  Let $m \in M$.  Then, we may write $m = er+fs+n$, where $r,s \in R$ and $n \in N$.  Again, we use Lemma \ref{translate} and apply $\phi$.
\begin{align}
\phi(umu^t) &= \begin{pmatrix} \alpha & \beta \\ \gamma & \delta \end{pmatrix} \begin{pmatrix} n & r \\ s & -n \end{pmatrix} \begin{pmatrix} \bar{\delta} & \bar{\beta} \\ \bar{\gamma} & \bar{\alpha} \end{pmatrix} \nonumber \\
&= \begin{pmatrix} \alpha n \bar{\delta} - \beta n \bar{\gamma} + s \beta \bar{\delta} + r \alpha \bar{\gamma} & \alpha n \bar{\beta} - \beta n \bar{\alpha} + s \beta \bar{\beta} + r \alpha \bar{\alpha} \\ \gamma n \bar{\delta} - \delta n \bar{\gamma} + s \delta \bar{\delta} + r \gamma \bar{\gamma} & \gamma n \bar{\beta} - \delta n \bar{\alpha} + s \delta \bar{\beta} + r \gamma \bar{\alpha} \end{pmatrix} \label{umu^t}
\end{align}

Recalling that $\alpha, \beta, \gamma, \delta \in Cl(N,q)$, the above is in $M$ for all $r,s \in R$ and $n \in N$ if and only if the remaining conditions (3)(iii) -- (vi) in the definition of $GV(N,q)$ hold.
\end{pf}

This also establishes the following:

\newtheorem{rep1}{Proposition \ref{vahlengroup}}
\begin{rep1}
Both $GV(N,q)$ and $SV(N,q)$ are indeed groups with inverse
\[
\begin{pmatrix} \alpha & \beta \\ \gamma & \delta \end{pmatrix}^{-1} = \frac{1}{\alpha \delta^t - \beta \gamma^t} \begin{pmatrix} \delta^t & -\beta^t \\ -\gamma^t & \alpha^t \end{pmatrix}.
\]
\qed
\end{rep1}

Using $\phi$ we may define a norm $N$ on $GV(N,q)$ by lifting the norm from $Cl(M,q)$.  We then observe that $N: GV(N,q) \to R^*$ is a homomorphism given by
\[
\begin{pmatrix} \alpha & \beta \\ \gamma & \delta \end{pmatrix} \mapsto \alpha \delta^t - \beta \gamma^t
\]
to which $SV(N,q)$ is the kernel.

The map $\phi$ also induces a grading on ${\rm MAT}_2(Cl(N,q))$ which agrees with the natural one on the matrices.

\begin{lem}
There is a $\mathbb{Z}_2$-grading on ${\rm MAT}_2(Cl(N,q))$, which agrees with the one induced from $M$ by $\phi$.  Even elements have the form $\begin{pmatrix} + & - \\ - & + \end{pmatrix}$ and odd elements $\begin{pmatrix} - & + \\ + & - \end{pmatrix}$, where $+$ and $-$ denote whether the entry is even, or odd, respectively.
\end{lem}

As usual, we denote the even and odd parts by using subscripts, e.g. we write ${\rm MAT}_2(Cl(N))_0$ and ${\rm MAT}_2(Cl(N))_1$ for the even and odd parts of ${\rm MAT}_2(Cl(N))$.  We have the following easy corollary.

\begin{cor}
Let $M$ be a module over $R$ with a quadratic form $q$.  Suppose that there exists a splitting $M = \la e,f \ra \perp N$ such that $\la e,f \ra$ is a hyperbolic line and $N$ is a free module.  Then,
\begin{align*}
GV(N,q)_0 &\cong NC_0(M,q) \\
SV(N,q)_0 &\cong Spin(M,q)
\end{align*}
\end{cor}

\section{Equivalent definitions of Vahlen groups}\label{sec:vahlenequiv}

We will give three different definitions of Vahlen groups and show that they are equivalent under certain conditions.  The first two definitions are generalisations of the classical definitions, while the third is the one that is given in Section \ref{sec:vahlen}.

The two new definitions that we will give contain a description of the type of elements that the entries of the matrices can be.  As was noted in the introduction, the entries are not necessarily invertible in general.  Indeed Maks gives an example over a real Clifford algebra of a Vahlen matrix which has no invertible entries \cite[p. 41]{maks}.   However, the set of possible entries do form a monoid.  Cnops in \cite{cnops} gave a characterisation of the monoid of entries in a (paravector) Vahlen group over a real Clifford algebra.

For our general setup of an arbitrary quadratic form $q$ over a commutative ring $R$, we describe a monoid in a different way:
\[
\T(M,q):=\{x \in Cl(M,q)-\{0\} : N(x) \in R \mbox{ and } xMx^t \subseteq M \}
\]
It is easy to see that the set $\T = \T(M,q)$ is a monoid; that is, it is closed, associative and contains an identity.  We note that the second condition is equivalent to $x'M \bar{x} \subseteq M$ and so $\T$ is closed under the grade automorphism.   If $u$ is in the Clifford group $CL(M,q)$, then clearly it is in $\T$ if and only if $N(u) \in R$.  In particular, if $\T$ is in fact a group, then it is a subgroup of the Clifford group.

%\begin{lem}\label{Nintdom}
%Let $R$ be an integral domain. If $N(\alpha), N(\bar{\alpha}) \in R$, then $N(\bar{\alpha}) =N(\alpha)$.
%\end{lem}
%\begin{pf}
%We have $N(\bar{\alpha})^2 = \bar{\alpha}\alpha\bar{\alpha}\alpha = \alpha\bar{\alpha}\bar{\alpha}\alpha = N(\alpha)N(\bar{\alpha})$.  Since $R$ is an integral domain, $N(\bar{\alpha}) =N(\alpha)$.
%\end{pf}

We may now prove the equivalence of three different definitions under some mild assumptions on $\T(M,q)$ and $R$.

\begin{thm}\label{vahlenequiv}
Let $N$ be a free module over a commutative integral domain $R$ with a quadratic form $q$.  Consider matrices
\[
\begin{pmatrix} \alpha & \beta \\ \gamma & \delta \end{pmatrix}  \in {\rm MAT}_2(Cl(N,q)).
\]
If $\T(N,q)$ is closed under the transpose map, then the following three sets of conditions are equivalent.
\begin{enumerate}[label = {\rm (\arabic*)}]
\item \begin{enumerate}[label = {\rm (\roman*)}]
	\item $\alpha, \beta, \gamma, \delta \in \T(N,q) \cup \{0\}$
	\item $\alpha \delta^t - \beta \gamma^t \in R^*$
	\item $\alpha\beta^t, \delta \gamma^t \in N \cup \{0\}$
	\end{enumerate}
\item \begin{enumerate}[label = {\rm (\roman*)}]
	\item $\alpha, \beta, \gamma, \delta \in \T(N,q) \cup \{0\}$
	\item $\alpha \delta^t - \beta \gamma^t \in R^*$
	\item $\bar{\alpha}\beta, \bar{\delta}\gamma \in N \cup \{0\}$
	\end{enumerate}
\item \begin{enumerate}[label ={\rm (\roman*)}]
	\item $\alpha \delta^t - \beta \gamma^t \in R^*$
	\item $\alpha \beta^t = \beta \alpha^t$ and $\gamma \delta^t = \delta \gamma^t$
	\item $\alpha \bar{\alpha}, \beta\bar{\beta}, \gamma \bar{\gamma}, \delta \bar{\delta} \in R$
	\item $\alpha\bar{\gamma}, \beta\bar{\delta} \in N \cup \{0\}$
	\item for all $x \in N$, $\alpha x \bar{\beta} + \beta \bar{x} \bar{\alpha}, \gamma x \bar{\delta} + \delta \bar{x} \bar{\gamma} \in R$
	\item for all $x \in N$, $\alpha x \bar{\delta} + \beta \bar{x} \bar{\gamma} \in N \cup \{0\}$
\end{enumerate}
\end{enumerate}
\end{thm}

Recall that we defined $GV(N,q)$ to be the third of the definitions above and we have already shown that it is in fact a group.  However, here we will only need that it is closed under taking inverses.  We assume for the rest of this section that $R$ is a commutative integral domain.  The proof of the above theorem will proceed with a series of lemmas.

\begin{lem}\label{23equiv}
If $\T(N,q)$ is closed under the transpose map, definitions {\rm (1)} and {\rm (2)} are equivalent.
\end{lem}
\begin{pf}
Consider the two equations:
\begin{align*}
\alpha(\bar{\alpha}\beta)\alpha^t &= N(\alpha) \beta\alpha^t \\
\beta^t(\beta'\bar{\alpha})\beta &= N(\beta^t) \bar{\alpha}\beta
\end{align*}
For the first, since $\alpha \in \T(N,q)$, we have that $\bar{\alpha}\beta \in M$ implies that $\beta\alpha^t \in M$.  Since $\T(N,q)$ is closed under the transpose map, $\beta^t \in \T(N,q)$ and we have the reverse implication from the second equation.
\end{pf}

The following lemma is an easy generalisation of a well-known result.

\begin{lem}{\rm \cite[Lemma 3.3]{vahlenpodd}}\label{link}
Let $g = \begin{pmatrix} \alpha & \beta \\ \gamma & \delta \end{pmatrix}$ such that $\alpha\delta^t - \beta \gamma^t \in R^*$ and $\alpha \bar{\alpha}, \beta\bar{\beta}, \gamma \bar{\gamma}, \delta \bar{\delta} \in R$.  If $\bar{\alpha}\beta, \bar{\gamma}\delta \in N$, then $\alpha^t\gamma, \beta^t\delta \in N$.
\end{lem}
\begin{pf}
The pseudo-determinant of $g$ is $r = \alpha\delta^t - \beta \gamma^t = \bar{\delta^t}\bar{\alpha} - \bar{\gamma^t}\bar{\beta}$.  Pre-multiply by $\gamma^t$ and post-multiply by $\alpha$ to obtain:
\[
r (\alpha^t \gamma)^t = r\gamma^t\alpha = (\bar{\delta}\gamma)^t(\bar{\alpha}\alpha) - (\bar{\gamma}\gamma)^t\bar{\beta}\alpha
\]
The result for $\beta^t\delta$ is similar.
\end{pf}

In particular, the above lemma is true for elements satisfying definition {\rm (2)}.  We will use this in the following lemma.

\begin{lem}\label{closedinv}
Let $\T(N,q)$ be closed under the transpose map.  Then, in definitions {\rm (1)} and {\rm (2)}, every element is invertible and its inverse, given by
\[
\begin{pmatrix} \alpha & \beta \\ \gamma & \delta \end{pmatrix}^{-1} = \frac{1}{\alpha \delta^t - \beta \gamma^t} \begin{pmatrix} \delta^t & -\beta^t \\ -\gamma^t & \alpha^t \end{pmatrix},
\]
satisfies the relations.  That is, the set is closed under taking inverses.
\end{lem}
\begin{pf}
By Lemma \ref{23equiv}, we can use either definition.  It is an easy calculation to see that the above is a right inverse for $g$ under definition (1).  To show it is also a left inverse consider the following.
\[
\begin{pmatrix} \delta^t & -\beta^t \\ -\gamma^t & \alpha^t \end{pmatrix}\begin{pmatrix} \alpha & \beta \\ \gamma & \delta \end{pmatrix} = \begin{pmatrix} \delta^t\alpha -\beta^t\gamma& \delta^t\beta - \beta^t\delta \\ \alpha^t\gamma -\gamma^t\alpha& \alpha^t\delta -\gamma^t\beta\end{pmatrix}
\]
From Lemma \ref{link}, the off diagonals are zero; this is also condition (iii) for the inverse.  Now observe that $(\alpha^t\delta -\gamma^t\beta)^t = \delta^t\alpha -\beta^t\gamma$, so it suffices to show that $\delta^t\alpha -\beta^t\gamma = \alpha\delta^t - \beta\gamma^t$; this being in $R^*$ is condition (i) for the inverse.  Using condition (iii) and Lemma \ref{link}, we have
\begin{align*}
\bar{\alpha}\alpha(\delta^t\alpha -\beta^t\gamma) &= \bar{\alpha}(\alpha\delta^t)\alpha -\bar{\alpha}\alpha\beta^t\gamma\\
&=\bar{\alpha}(\alpha\delta^t)\alpha -\bar{\alpha}\beta\alpha^t\gamma\\
&=\bar{\alpha}(\alpha\delta^t)\alpha -\bar{\alpha}\beta\gamma^t\alpha\\
&=\bar{\alpha}(\alpha\delta^t -\beta\gamma^t)\alpha\\
&=\bar{\alpha}\alpha(\alpha\delta^t -\beta\gamma^t)
\end{align*}
Since $\bar{\alpha} \in \T(N,q)$, $N(\bar{\alpha}) = \bar{\alpha}\alpha \in R$.  So, we see that $(\delta^t\alpha -\beta^t\gamma)(\alpha\delta^t -\beta\gamma^t)^{-1} = 1$.  Hence, $g^{-1}$ is indeed also a left inverse and it satisfies condition (ii).  By assumption, $\T(N,q)$ is closed under the transpose map, so $g^{-1}$ also satisfies condition (i) and the proof is complete.
\end{pf}

%Let $\mu$ denote the action of $g$ on the left.  By the existence of a right inverse, there is a $\nu$ such that $\mu \circ \nu = id$ and so $\mu$ is surjective.  Now, $R$ is commutative and hence Noetherian.  So, $A:=Cl(N,q)$, which is finitely generated as an algebra, is also Noetherian.  By considering the chain of submodules $\ker \mu \subseteq \ker \mu^2 \subseteq \dots$, we see that $\ker \mu^k = \ker \mu^{k+1}$ for some $k$.  Hence, $\mu(\mu^ka) = 0$ implies that $\mu^ka = 0$, for $a \in A$.  That is, $\mu$ restricted to $\mu^kA$ is injective.  However, $\mu$ is surjective, so $\mu^kA = A$ and hence $\mu$ is an isomorphism.  Therefore, $g$ has a well-defined two-sided inverse.

%It remains to show that the inverse satisfies the relations.  By Lemma \ref{Cliffclosed}, $\alpha^t, \beta^t, \gamma^t, \delta^t \in C(N,q)$.  Calculating $g^{-1}g$, we see that this equals the identity only if $\delta^t\alpha - \beta^t\gamma \in R^*$.  However, this is condition (ii) for the inverse.  Lemma \ref{link} used on definition (2) gives $\alpha^t\gamma, \beta^t\delta^t \in N \cup \{0\}$, which is condition (1)(iii).

\begin{lem}
Suppose that $g:=\begin{pmatrix} \alpha & \beta \\ \gamma & \delta \end{pmatrix}$ satisfies either definition {\rm (3)}, or satisfies definitions {\rm (1)}, or {\rm (2)}, and $\T(N,q)$ is closed under the transpose map.  Then, $N(\alpha)=N(\bar{\alpha})$, $N(\beta)=N(\bar{\beta})$, $N(\gamma)=N(\bar{\gamma})$ and $N(\delta)=N(\bar{\delta})$.
\end{lem}
\begin{pf}
We will show the result for $\alpha$; the others follow similarly.  We may suppose that $\alpha \neq 0$ which using the pseudo-determinant condition implies that $\beta,\gamma \neq 0$.  In all cases, $g$ is invertible, so the pseudo-determinant condition for $g$ and $g^{-1}$ give us $\alpha\delta^t - \beta\gamma^t = \delta^t\alpha - \beta^t \gamma$.   Multiply by $\bar{\alpha}$ to obtain
\[
\bar{\alpha}(\alpha\delta^t - \beta\gamma^t) = (\delta^t\alpha - \beta^t \gamma)\bar{\alpha}
\]
Rearranging, we get
\begin{align*}
(\bar{\alpha}\alpha - \alpha \bar{\alpha})\delta^t &= \bar{\alpha}\beta\gamma^t - \beta^t\gamma\bar{\alpha}\\
&= \beta^t\bar{\alpha}^t\gamma^t - \beta^t\gamma\bar{\alpha}\\
&= \beta^t\left( (\gamma\bar{\alpha})^t - \gamma\bar{\alpha} \right)\\
&= \beta^t\left( \gamma\bar{\alpha} - \gamma\bar{\alpha} \right) = 0
\end{align*}

\end{pf}

We can now complete the proof.  This is the only place where we need to use that definition (3) is closed under taking inverses.

\begin{pf}[Proof of Theorem $\ref{vahlenequiv}$]
It remains to show that (1), or (2), are equivalent to (3).  Clearly, definition (1) implies conditions (i)--(iii) of definition (3).  Condition (3)(iv) is just (2)(iii) applied to the inverse.  Assume definition (1) and suppose $\alpha, \beta \neq 0$.  Then,
\[
N(\beta^t)\alpha n \bar{\beta} = \alpha(\beta^t\beta')n\bar{\beta} = \alpha\beta^t(\beta n' \beta^t)'
\]
In definition (1), $\beta \in \T(N,q)$ and $\alpha\beta^t \in N \cup \{0\}$, so $x:=\alpha n \bar{\beta}$ is the product of two elements of $N$.  Hence, by the Clifford relations, $\alpha n \bar{\beta} + \beta \bar{n} \bar{\alpha} = x + \bar{x} = x + x^t\in R$.  Similarly, $\gamma n \bar{\delta} + \delta \bar{n} \bar{\gamma} \in R$ and condition (3)(v) is satisfied.

We have shown that conditions (3)(i)--(v) are satisfied by both definitions (1) and (3).  So under these assumptions, it is enough to show the equivalence of (1)(i) and (3)(vi).  By (3)(iii), we need just show that $\alpha n \alpha^t \in N$ is equivalent to $\alpha n \bar{\delta} - \beta n \bar{\gamma} \in N$ for $\alpha \neq 0$.

%We begin by claiming that $N(\alpha)\bar{\delta} = N(\alpha^t)\bar{\delta}$.  All definitions are closed under inverses, so the pseudo-determinant condition for $g$ and $g^{-1}$ give us $\alpha\delta^t - \beta\gamma^t = \delta^t\alpha - \beta^t \gamma$.   Multiply by $\bar{\alpha}$ to obtain
%\[
%\bar{\alpha}(\alpha\delta^t - \beta\gamma^t) = (\delta^t\alpha - \beta^t \gamma)\bar{\alpha}
%\]
%Rearranging, we get
%\begin{align*}
%(\bar{\alpha}\alpha - \alpha \bar{\alpha})\delta^t &= \bar{\alpha}\beta\gamma^t - \beta^t\gamma\bar{\alpha}\\
%&= \beta^t\bar{\alpha}^t\gamma^t - \beta^t\gamma\bar{\alpha}\\
%&= \beta^t\left( (\gamma\bar{\alpha})^t - \gamma\bar{\alpha} \right) = 0
%\end{align*}
%Since $N(\alpha), N(\alpha^t) \in R$, applying the grade automorphism to both sides proves the claim.  We may now proceed with the remainder of the proof.
\begin{align*}
N(\alpha^t)(\alpha n \bar{\delta} + \beta \bar{n} \bar{\gamma}) &= \alpha n N(\alpha^t)\bar{\delta} + N(\alpha^t)\beta \bar{n}\bar{\gamma} \\
&=\alpha n \alpha^t(\alpha\delta^t - \beta\gamma^t)' + \alpha n \alpha^t\beta'\bar{\gamma} + N(\alpha^t)\beta \bar{n} \bar{\gamma}
\end{align*}
If $\beta = 0$, then our proof is complete.  So, we may assume that $\beta \neq 0$.  Now, $\alpha^t\beta' = (\bar{\alpha}\beta)' \in N$.  Since $n$ is also in $N$, by the Clifford relations, $n\alpha^t\beta' + \alpha^t\beta'n = r \in R$.  So, $n\alpha^t\beta' =  r +  \alpha^t\beta'n = r - \overline{\alpha^t\beta'}\bar{n} = r - \beta^t\alpha'\bar{n}$.  We have
\begin{align*}
N(\alpha^t)(\alpha n \bar{\delta} + \beta \bar{n} \bar{\gamma}) &= \alpha n \alpha^t(\alpha\delta^t - \beta\gamma^t)' + \alpha (r - \beta^t\alpha' \bar{n})\bar{\gamma} + N(\alpha^t)\beta \bar{n} \bar{\gamma} \\
&= \alpha n \alpha^t(\alpha\delta^t - \beta\gamma^t)' + r \alpha\bar{\gamma} - \beta\alpha^t\alpha' \bar{n}\bar{\gamma} + N(\alpha^t)\beta \bar{n} \bar{\gamma} \\
& =  \alpha n \alpha^t(\alpha\delta^t - \beta\gamma^t)' + r\alpha\bar{\gamma}
\end{align*}
Now, $\alpha\bar{\gamma} \in N \cup \{0\}$ in definition (3) and also in definition (2)(iii) applied to the inverse.   so we have proved our claim.  The analogous results for $\beta$, $\gamma$ and $\delta$ follow in a similar way.
\end{pf}

We note that when $R$ is restricted to $\mathbb{R}$, or $\mathbb{C}$ and $q$ is a definite quadratic form, then the conditions are satisfied and definition {\rm (1)}, or {\rm (2)}, of $GV(N,q)$ are the two usual definitions.

\section{Paravector Vahlen groups}\label{sec:paravahlen}

Let $L$ be a free module over a commutative ring $R$ with a quadratic form $q$.  Then, the Clifford algebra $Cl(L,q)$ is defined.  A \emph{paravector} is an element of $R \oplus L$.

\begin{defn}
Let $L$ be a free module over a commutative ring $R$ with a quadratic form $q$.  We define
\[
GPV(L,q):=\left\{ \begin{pmatrix} \alpha & \beta \\ \gamma & \delta \end{pmatrix}  \in {\rm MAT}_2(Cl(L,q)) \right\}
\]
such that
\begin{enumerate}[label = (\roman*)]
\item $\alpha \delta^t - \beta \gamma^t \in R^*$
\item $\alpha \beta^t = \beta \alpha^t$ and $\gamma \delta^t = \delta \gamma^t$
\item $\alpha \bar{\alpha}, \beta\bar{\beta}, \gamma \bar{\gamma}, \delta \bar{\delta} \in R$
\item $\alpha\bar{\gamma}, \beta\bar{\delta} \in R\oplus L$
\item for all $x \in R \oplus L$, $\alpha x \bar{\beta} + \beta \bar{x} \bar{\alpha}, \gamma x \bar{\delta} + \delta \bar{x} \bar{\gamma} \in R$
\item for all $x \in R \oplus L$, $\alpha x \bar{\delta} + \beta \bar{x} \bar{\gamma} \in R \oplus L$
\end{enumerate}
We call $\alpha \delta^t - \beta \gamma^t$ the \emph{pseudo-determinant} and define $SPV(L,q)$ to be the subset of $GPV(L,q)$ with pseudo-determinant $1$.
\end{defn}

We will call the above groups the \emph{general paravector Vahlen group} and \emph{special paravector Vahlen group}, respectively.  However, when the context is clear, we will just refer to them as Vahlen groups.  We note that this extends the definition of Vahlen groups given by Elstrodt, Grunewald and Mennicke in \cite{vahlen} from fields in odd characteristic $p$ to arbitrary commutative rings.

As before, we will show that $GPV(L,q)$ and $SPV(L,q)$ are indeed groups.

\begin{prop}\label{GPVgroup}
The sets $GPV(L,q)$ and $SPV(L,q)$ are groups with inverse
\[
\begin{pmatrix} \alpha & \beta \\ \gamma & \delta \end{pmatrix}^{-1} = \frac{1}{\alpha \delta^t - \beta \gamma^t} \begin{pmatrix} \delta^t & -\beta^t \\ -\gamma^t & \alpha^t \end{pmatrix}.
\]
Furthermore, the pseudo-determinant is a group homomorphism into $R^*$, whose kernel is $SPV(L,q)$.
\end{prop}

We will proceed in a similar way as for the analogous proposition for $GV(N,q)$.  First we recall a standard construction which can be found, for example, in \cite[Chapter 7]{omeara}. 

Suppose that $M = \la z \ra \perp Z$ is a splitting of $M$, where $Z$ is a free submodule and $q(z) = s \in R^*$.  Then, we may define a map $\psi: Z \to Cl(M,q)$ by $x \mapsto zx$.  Then,
\[
\psi(x)^2 = zxzx = -sx^2 = -s q(x)
\]
Let $Z^{-s}$ denote the module $Z$ where the quadratic form $q_z$ has been adjusted by scaling by $-s$.  Then, $(Cl_0(M,q),\psi)$ is compatible with $Z^{-s}$.  Hence, this defines an algebra homomorphism $\psi_z:Cl(Z^{-s}, q_z) \to Cl_0(M,q)$.

Recall our previous assumptions on $M$, that $M = \la e,f \ra \perp N$.  We now assume further that $M = \la e,f \ra \perp \la z \ra \perp L$, where $L$ is a free submodule, $\la e,f \ra$ is a hyperbolic line and $q(z) = -1\in R^*$.  Define $Z = \la e,f \ra \perp L$.  Noting that $s=-1$, $Z^{-s} = Z$ and $q_z = q$.  From before, we have an isomorphism $\phi:Cl(Z,q) \to {\rm MAT}_2(Cl(L,q))$.  So, writing $\psi = \psi_z$, we may now define $\theta = \phi \circ \psi^{-1}: Cl_0(M,q) \to {\rm MAT}_2(Cl(L,q))$.

\begin{lem}
The isomorphism $\theta = \phi \circ \psi^{-1}: Cl_0(M,q) \to {\rm MAT}_2(Cl(L,q))$ is given by
\[
\psi(\alpha) ef + \psi(\beta) ez + \psi(\gamma')fz + \psi(\delta')fe \mapsto \begin{pmatrix} \alpha & \beta \\ \gamma & \delta \end{pmatrix}
\]
\end{lem}
\begin{pf}
It is an easy calculation to write down $\theta^{-1}$ which gives the above mapping.
\end{pf}

Before continuing, we first observe a simple lemma about $\psi$.

\begin{lem}
Let $y \in Cl(Z,q)$.  We may write $y = y_0 + y_1$ where $y_0 \in Cl(Z,q)_0$ and $y_1 \in Cl(Z,q)_1$.  Then,
\[
\psi(y) = y_0 + zy_1
\]
\end{lem}
\begin{pf}
Observe that, if $u,v \in Z$, then $\psi(uv) = uzvz = -uz^2v = -q(z)uv = uv$.  Hence, $\psi$ acts as the identity on even products.
\end{pf}

\begin{lem}\label{psibar}
Let $y \in Cl(Z,q)$.  Then,
\[
\overline{\psi(y)} = \psi(\bar{y})
\]
\end{lem}
\begin{pf}
Let $y = y_0 + y_1$ be the graded decomposition.  Then,
\[
\overline{\psi(y)} =\psi(y)^t =  (y_0 + zy_1)^t = y_0^t+y_1^tz = y_0^t - zy_1^t = \psi(\bar{y}).
\]
\end{pf}

\begin{lem}\label{GVbar}
Suppose that $u \in Cl_0(M,q)$ such that $\theta(u) = \begin{pmatrix} \alpha & \beta \\ \gamma & \delta \end{pmatrix}$.  Then,
\[
\theta(\bar{u}) =\theta(u^t)= \begin{pmatrix} \delta^t & -\beta^t \\ -\gamma^t & \alpha^t \end{pmatrix}
\]
\end{lem}
\begin{pf}
By Lemmas \ref{psibar} and \ref{translate}, $\theta(\bar{u}) = \phi(\psi^{-1}(\bar{u})) = \phi(\overline{\psi^{-1}(u)}) = \overline{\phi(\psi^{-1}(u))}$.  Hence,
\[
\theta(\bar{u}) =\overline{\theta(u)}= \begin{pmatrix} \delta^t & -\beta^t \\ -\gamma^t & \alpha^t \end{pmatrix}
\]
\end{pf}

We may now proceed with the main theorem for this section.  Recall that we defined
\[
NC(M,q):=\{ u \in C(M,q) : N(u) \in R^* \}
\]
and so $NC_0(M,q)$ is the even part of $NC(M,q)$.

\begin{thm}\label{paravahleniso}
Let $M$ be a module over a commutative ring $R$ with a quadratic form $q$.  Suppose that there exists a splitting $M = \la e,f \ra \perp \la z \ra \perp L$, where $L$ is a free submodule, $\la e,f \ra$ is a hyperbolic line and $q(z) = -1$.  Then,
\begin{align*}
GPV(L,q) &\cong GV(L \perp \la z \ra,q)_0 \cong NC_0(M,q)\\
SPV(L,q) &\cong SV(L \perp \la z \ra,q)_0 \cong Spin(M,q)
\end{align*}
\end{thm}
\begin{pf}
The proof will run analogously to that of Theorem \ref{vahleniso}, indeed some of the calculations are identical and we will borrow them from there.  To start, if $u \in NC_0(M,q)$, then $N(u) \in R^*$ and the inverse of $u$ is $\bar{u}/N(u)$.   So again we calculate $\theta(u\bar{u})$.  By Lemma \ref{GVbar}, this is exactly the same calculation as in equation (\ref{ubaru}) and so we get conditions (i) and (ii) of the definition of $GPV(L,q)$.

Now, as $u$ is even, $u \in C_0(M,q)$ if and only if $u'm\bar{u} = um\bar{u} \in M$ for all $m \in M$.  Since $z$ is invertible and hence has a trivial annihilator,
\[
 um\bar{u} \in M \quad \iff \quad zum\bar{u} \in zM
\]
Recalling that $M = \la z \ra \perp Z$, we see that $zM = z\la z \ra \perp zZ = R \oplus zZ$, which has even grading.  In particular, it lies in the image of $\psi$ and $\psi^{-1}(zM) = R \oplus Z$.  Fixing notation, given $m \in M$, let $y \in R \oplus Z$ such that $\psi(y) = zm$.  Define $x:=\phi^{-1}(u) \in Cl(Z,q)$.  If $x$ has graded decomposition $x=x_0 + x_1$, then $u = \phi(x) = x_0 + zx_1$.  Now, we have
\[
zum\bar{u} = z(x_0+zx_1)m\bar{u} = (x_0 - zx_1)(zm)\bar{u} = \psi(x')\psi(y)\psi(\bar{x}) = \psi(x'y\bar{x})
\]
Hence, $u \in C_0(M,q)$ if and only if $x'y\bar{x} \in R \oplus Z$ for all $y \in R \oplus Z$ which is equivalent to $xyx^t \in R \oplus Z$.  We also performed this calculation before in equation (\ref{umu^t}) of the proof of Theorem \ref{vahleniso}.  However, we must now allow $n \in R \oplus L$ and alter the conclusion to allow the diagonals to be in $R \oplus L$.  Hence, we obtain condition (iii)--(vi) for the paravector Vahlen groups.
\end{pf}

This also establishes the proof of Proposition \ref{GPVgroup}, showing that both $GPV(L,q)$ and $SPV(L,q)$ are indeed groups.

\section{Equivalent definitions of paravector Vahlen groups}\label{sec:paravahlenequiv}

We will now give two other equivalent definitions of paravector Vahlen groups.  But before we do, we define
\[
\mathcal{PT}(L,q) := \{ x \in C(L,q) : N(x) \in R \mbox{ and } x(R \oplus L)x^t \subseteq R \oplus L \}
\]
As before, this is a monoid which is closed under the grade automorphism.  We note in passing that if $u$ is in the so-called paravector Clifford group (see \cite[Chapter 13]{porteoustopgeom} for details), then it is in $\mathcal{PT}(L,q)$ if and only if $N(u) \in R$.  If $\mathcal{PT}(L,q)$ is a group, then it is a subgroup of the paravector Clifford group.

\begin{thm}\label{paravahlenequiv}
Let $L$ be a free module over a commutative integral domain $R$ with a quadratic form $q$.  Consider matrices
\[
\begin{pmatrix} \alpha & \beta \\ \gamma & \delta \end{pmatrix}  \in {\rm MAT}_2(Cl(L,q)).
\]
If $\mathcal{PT}(L,q)$ is closed under the transpose map, then the following three sets of conditions are equivalent.
\begin{enumerate}[label = {\rm (\arabic*)}]
\item \begin{enumerate}[label = {\rm (\roman*)}]
	\item $\alpha, \beta, \gamma, \delta \in \mathcal{PT}(L,q) \cup \{0\}$
	\item $\alpha \delta^t - \beta \gamma^t \in R^*$
	\item $\alpha\beta^t, \delta \gamma^t \in R \oplus L$
	\end{enumerate}
\item \begin{enumerate}[label = {\rm (\roman*)}]
	\item $\alpha, \beta, \gamma, \delta \in \mathcal{PT}(L,q) \cup \{0\}$
	\item $\alpha \delta^t - \beta \gamma^t \in R^*$
	\item $\bar{\alpha}\beta, \bar{\delta}\gamma \in R \oplus L$
	\end{enumerate}
\item \begin{enumerate}[label = {\rm (\roman*)}]
	\item $\alpha \delta^t - \beta \gamma^t \in R^*$
	\item $\alpha \beta^t = \beta \alpha^t$ and $\gamma \delta^t = \delta \gamma^t$
	\item $\alpha \bar{\alpha}, \beta\bar{\beta}, \gamma \bar{\gamma}, \delta \bar{\delta} \in R$
	\item $\alpha\bar{\gamma}, \beta\bar{\delta} \in R\oplus L$
	\item for all $x \in R \oplus L$, $\alpha x \bar{\beta} + \beta \bar{x} \bar{\alpha}, \gamma x \bar{\delta} + \delta \bar{x} \bar{\gamma} \in R$
	\item for all $x \in R \oplus L$, $\alpha x \bar{\delta} + \beta \bar{x} \bar{\gamma} \in R \oplus L$
\end{enumerate}
\end{enumerate}
\end{thm}

As before, the third definition is what we defined $GPV(L,q)$ to be, which we have already shown to be a group.  Assume for the rest of this section that $R$ is a commutative integral domain.  The proof of the above theorem will proceed in the same way as that for Theorem \ref{vahlenequiv}, via a series of lemmas.

\begin{lem}\label{para23equiv}
Definitions {\rm (1)} and {\rm (2)} are equivalent.
\end{lem}
\begin{pf}
The proof is the same as for Lemma \ref{23equiv}.
\end{pf}

\begin{lem}{\rm \cite[Lemma 3.3]{vahlenpodd}}\label{paralink}
Let $g = \begin{pmatrix} \alpha & \beta \\ \gamma & \delta \end{pmatrix}$ such that $\alpha\delta^t - \beta \gamma^t \in R^*$ and $\alpha \bar{\alpha}, \beta\bar{\beta}, \gamma \bar{\gamma}, \delta \bar{\delta} \in R$.  If $\bar{\alpha}\beta, \bar{\gamma}\delta \in R \oplus L$, then $\alpha^t\gamma, \beta^t\delta \in R \oplus L$.
\end{lem}
\begin{pf}
The proof is the same as for Lemma \ref{link}
\end{pf}

\begin{lem}
In definitions {\rm (1)} and {\rm (2)}, every element is invertible and its inverse, given by
\[
\begin{pmatrix} \alpha & \beta \\ \gamma & \delta \end{pmatrix}^{-1} = \frac{1}{\alpha \delta^t - \beta \gamma^t} \begin{pmatrix} \delta^t & -\beta^t \\ -\gamma^t & \alpha^t \end{pmatrix},
\]
satisfies the relations.  That is, the set is closed under taking inverses.
\end{lem}
\begin{pf}
The proof is the same as for Lemma \ref{closedinv}
\end{pf}

We can now complete the proof.  This is the only place where we need to use that definition (3) is closed under taking inverses.

\begin{pf}[Proof of Theorem $\ref{paravahlenequiv}$]
This proof runs in a similar way to that of Theorem \ref{vahlenequiv}, but with some changes due to the paravectors.  It remains to show that (1), or (2), are equivalent to (3).  As before, definitions (1) and (2) implies conditions (i)--(iv) of definition (3).  Assume definition (1) and suppose $\alpha, \beta \neq 0$, then as before
\[
N(\beta^t)\alpha x \bar{\beta} = \alpha(\beta^t\beta')x\bar{\beta} = \alpha\beta^t(\beta x' \beta^t)'
\]
So $y:=\alpha x \bar{\beta}$ is the product of two elements of $R \oplus L$.  Such a product lies in $R \oplus L \oplus (L\otimes L)$.  Let $y = r + l + nm \in R \oplus L \oplus (L\otimes L)$.  Then
\[
y + \bar{y} = r + l + mn + r -l + nm = 2r + mn+nm \in R
\]
Hence, $\alpha x \bar{\beta} + \beta \bar{x} \bar{\alpha} = y + \bar{y} \in R$.  Similarly, $\gamma x \bar{\delta} + \delta \bar{x} \bar{\gamma} \in K$ and condition (3)(v) is satisfied.

We have shown that conditions (3)(i)--(v) are satisfied by both definitions (1) and (3).  So, as before, it is enough to show that $\alpha x \alpha^t \in N$ is equivalent to $\alpha x \bar{\delta} + \beta \bar{x} \bar{\gamma} \in N$ for $\alpha \neq 0$.
\begin{align*}
N(\alpha^t)(\alpha x \bar{\delta} + \beta \bar{x} \bar{\gamma}) &= \alpha x \alpha^t\alpha'\bar{\delta} + N(\alpha^t)\beta \bar{x} \bar{\gamma} \\
&=\alpha x \alpha^t(\alpha \delta^t - \beta\gamma^t)' + \alpha x \alpha^t\beta'\bar{\gamma} + N(\alpha^t)\beta \bar{x} \bar{\gamma}
\end{align*}
If $\beta = 0$, then our proof is complete.  So, we may assume that $\beta \neq 0$.  Observe that $\alpha^t\beta', x \in R \oplus L$.  So, by the above, $x(\alpha^t\beta') + \overline{(\alpha^t\beta')}\bar{x} =:r \in R$.  Hence, $x(\alpha^t\beta') = r - \beta^t\alpha'\bar{x}$.  We have
\begin{align*}
N(\alpha^t)(\alpha x \bar{\delta} + \beta \bar{x} \bar{\gamma}) &= \alpha x \alpha^t(\alpha\delta^t - \beta\gamma^t)' + \alpha (r - \beta^t\alpha' \bar{x})\bar{\gamma} + N(\alpha^t)\beta \bar{x} \bar{\gamma} \\
& =  \alpha x \alpha^t(\alpha\delta^t - \beta\gamma^t)' + r\alpha\bar{\gamma}
\end{align*}
Now, $\alpha\bar{\gamma} \in N \cup \{0\}$ in definition (3) and also in definition (2)(iii) applied to the inverse.   so we have proved our claim.  The analogous results for $\beta$, $\gamma$ and $\delta$ follow in a similar way.
\end{pf}

\end{document}